\documentclass[10pt,twoside]{article}
\usepackage[psamsfonts]{amssymb}
\usepackage{amsfonts}
\usepackage[mathcal]{eucal}
\usepackage{Latex-document}

\newcommand{\DD}[3]{\|Df^{#1}_{/E({#2})}\|\|Df^{-#1}_{/F(f^{#1}({#2}))}\|^{#3}}
\newcommand{\noi}{\noindent}
\newcommand{\Li}{\tilde{\cal L}}
\newcommand{\e}{\epsilon}
\newcommand{\de}{\delta}
\newcommand{\La}{\Lambda}

\title{\bf Tangent Bundles Dynamics and Its\vskip -2mm Consequences\thanks{Partially supported by FAPERJ-Brazil
and Guggenheim Foundation.}\vskip 6mm}
\author{E. R. Pujals\vspace*{-0.5cm}\thanks{Instituto de Matem\'atica, Universidade Federal do Rio de Janeiro,
C. P. 68.530, CEP 21.945-970, Rio de Janeiro, R. J. , Brazil.
E-mail: enrique@im.ufrj.br}}
\date{\vspace{-8mm}}

\begin{document}

\maketitle

\thispagestyle{first} \setcounter{page}{327}

\begin{abstract}

\vskip 3mm

We will consider here some dynamics of the tangent map, weaker
than hyperbolicity, and we will discuss if these structures are
rich enough to provide a good description of the dynamics from a
topological and geometrical point of view. This results are
useful in attempting to obtain global scenario in terms of
generic phenomena relative both to the space of dynamics and to
the space of trajectories. Moreover, we will relate these results
with the study of systems that remain globally transitive under
small perturbations. \vskip 4.5mm


\end{abstract}

\vskip 12mm

\section{Introduction}

\vskip-5mm \hspace{5mm}

A long time goal in the theory of dynamical systems is to describe
the dynamics of ``big sets'' (generic or residual, dense, etc) in
the space of all dynamical systems.

It was thought in the sixties that this could be realized by the
so called hyperbolic ones: systems with the assumption that the
tangent bundle over the Limit set ($L(f)$, the accumulation points
of any orbit) splits into two complementary subbundles that are
uniformly forward (respectively backward) contracted by the
tangent map by. The richness of this description would follow from
the fact that {\em the hyperbolic dynamic on the tangent bundle
characterizes the dynamic over the manifold} from a geometrical
and topological point of view.

Nevertheless, uniform hyperbolicity were soon realized to be a property less universal than it was initially
thought: {\em there are open sets in the space of dynamics which are non-hyperbolic.} After some initial examples
of non-density of the hyperbolic systems in the universe of all systems (see \cite{S,AS}), two key aspects were
focused in these examples. On one hand, open sets of non-hyperbolic diffeomorphisms which remain transitive under
perturbation (existence of a dense orbit for any system). On the other, residual sets of non-hyperbolic
diffeomorphisms, each one exhibiting infinitely many transitive sets. Roughly speaking, it was showed $\,$ that
two kind of different phenomena $\,$ can appear $\,$in the complement of the hyperbolic systems: {\bf a)} {\em
dynamics that robustly can be decomposed into a finite number of closed transitive sets;} {\bf b)} {\em dynamics
that generically exhibit infinitely many disjoint transitive sets.}

The first kind of phenomena occurs in dimension higher or equal
than $3$ and although the examples are not hyperbolic they exhibit {\em
some kind of decomposition of the tangent bundle into invariant
subbundles.} The second one, was obtained by Newhouse (see
\cite{N1}, \cite{N2}, \cite {N3}), who following an early work of
the non-density of hyperbolicity for $C^2$ surface maps, showed
that the unfolding of a {\em homoclinic bifurcation} (non
transversal intersection of stable and unstable manifolds of
periodic points) leads to a very rich dynamics: residual subsets
of open sets of diffeomorphisms whose elements display infinitely
many sinks.

These new results naturally pushed some aspects of the theory on
dynamical systems in different directions:

\begin{enumerate}
\item The study of the dynamical phenomena obtained from homoclinic
bifurcations;
\item  The characterization of universal mechanisms that could yield
to robustly non-hyperbolic behavior;
\item The study and characterization of isolated transitive
sets that remain transitive for all nearby system ({\it robust
transitivity});
\item The dynamical consequences that follows from
some kind of the dynamics over the tangent bundle, weaker than the
hyperbolic one.
\end{enumerate}

As we will show, these problems are related
and they indeed constitute different aspects of the same
phenomena. In many cases, such relations provide a conceptual
framework, as the hyperbolic theory did for the case of transverse
homoclinic orbits.

In the next section, we will discuss the previous aspects for the case of surfaces maps, and in particular we will
consider a dynamics in the tangent bundle weaker than hyperbolicity, called dominated splitting. In section
\ref{robust} we will discuss the problems about the robust transitivity and its relation with other dynamics on
the tangent bundle. Finally, in the last section, we will consider the equivalent problems for flows taking into
account their intrinsic characteristic that leads to further questions and difficulties do the presence of
singularities. We point out that in this survey we will focus more into topological and geometrical aspects of the
dynamic rather on the ergodic ones.

Many of the issues discussed here are consequences of works and
talks with Martin Sambarino. I also want to thanks Maria J.
Pacifico for her help to improve this article.

\section{Surfaces maps, homoclinic tangencies and ``non-critical"
behaviors} \label{surfaces}

\vskip-5mm \hspace{5mm}

After the seminal works of Newhouse,  many others were developed
in the direction to understand the phenomena that could appear
after a bifurcation of homoclinic tangencies (tangent intersection
of stable and unstable manifolds of periodic points). In fact,
other fundamental dynamic prototype were found in this context,
namely the so called cascade of bifurcations, the H\'enon-like
strange attractor (\cite{BC}, \cite{MV}) (even infinitely many
coexisting ones \cite{C}), and superexponential growth of periodic
points (\cite{K}). Even before these last results, Palis
(\cite{PT}, \cite{P}) conjectured that the presence of a
homoclinic tangency is a very common phenomenon in the complement
of the closure of the hyperbolic ones. In fact, if the conjecture
is true, then homoclinic bifurcation could play a central role in
the global understanding of the space of dynamics for it would
imply that each of these bifurcation phenomena is dense in the
complement of the closure of the hyperbolic ones. More precisely,
he conjectured that {\em Every $f\in Diff^r(M^2), r\ge 1,$ can be
$C^r$-approximated by a diffeomorphism exhibiting either a
homoclinic tangency or by one which is  hyperbolic.}

The presence of homoclinic tangencies have many analogies with the
presence of {\em critical points} for one-dimensional
endomorphisms. Homoclinic tangecies correspond in the one
dimensional setting to preperiodic critical points and it is known
that its bifurcation leads to complex dynamics. On the other hand,
Ma\~n\'e  (see \cite{M1})showed that for regular and generic
one-dimensional endomorphisms, the {\it absence of critical points
is enough to guarantee hyperbolicity.} This result raises the
question about the dynamical properties of surface maps exhibiting
no homoclinic tangencies. In this direction, first it is proved in
that some kind of dynamic over the tangent bundle (weaker than the
hyperbolic one) can be obtained in the robust lack of homoclinic
tangencies. And later, it is showed that this dynamic on the
tangent bundle is rich enough to describe the dynamic on the
manifold. More precisely:

{\bf Theorem 1} (\cite{PS1}){\bf: } {\em Surface diffeomorphisms
that can not be $C^1$-approximated by another exhibiting
homoclinic tangencies, has the property that its Limit set has}
{\bf dominated splitting}. \vskip 2pt An $f$-invariant set
$\Lambda$ has dominated splitting if the tangent bundle can be
decomposed into two invariant subbundles $T_\Lambda M=E\oplus F
,$ such that:
$$\DD{n}{x}{}\le C\lambda^n, \mbox{ for all } x\in\Lambda, n\ge 0,$$

\noi with $C>0$ and $0<\lambda<1.$

As, dominated splitting prevents the presence of tangencies, we
could say that domination plays for surface diffeomorphisms the
role that the {\it non-critical behavior} does for one dimensional
endomorphisms.

To have a satisfactory  description for this non-critical
behavior (existence of a dominated splitting),
we should describe its dynamical consequences. A natural question arises:
{\it is it possible to describe the dynamics of a
system having dominated splitting?}

The next result gives a positive answer (as satisfactory as in
hyperbolic case) when $M$ is a compact surface. More precisely, we
give a complete description of the topological dynamics of a
$C^2$ system having a dominated splitting.
Actually, first, the dominated decomposition is
understood under a generic assumption.

{\bf Theorem 2} (\cite{PS1}) {\bf : } {\em Let $f\in Diff^2(M^2)$
and assume that $\Lambda \subset L(f)$ is a compact invariant set
exhibiting a dominated splitting such that any periodic point is
a hyperbolic periodic point. Then,
$\Lambda=\Lambda_1\cup{\Lambda_2}$ where $\Lambda_1$ is hyperbolic
and $\Lambda_2$ consists of a finite union of periodic simple
closed curves ${\cal{C}}_1,...{\cal{C}}_n$, normally hyperbolic,
and such that $f^{m_i}:{\cal{C}}_i\to {\cal{C}}_i$ is conjugated
to an irrational rotation ($m_i$ denotes the period of
${\cal{C}}_i$).}

Using this Theorem and understanding the obstruction for the
hyperbolicity assuming domination, we can
characterize $L(f)$ without any generic assumption.

{\bf Theorem 3} (\cite{PS2}) {\bf :} {\em Let $f\in Diff^2(M^2)$ and
assume that $L(f)$ has a dominated splitting. Then $L(f)$ can be
decomposed into $L(f)={\mathcal{I}}\cup{\Li(f)}\cup{\mathcal{R}}$
such that:

$1 .$ $\mathcal{I}$ is contained in a finite union of normally
hyperbolic periodic arcs.

$2 .$ $\mathcal{R}$ is a finite union of normally hyperbolic periodic
simple closed curves supporting an irrational rotation.

$3 .$ $f{/\Li(f)}$ is expansive and admits a spectral
decomposition (into finitely many homoclinic classes).}

Roughly speaking, the above theorem says that the dynamics of a
$C^2$ diffeomorphism having a dominated splitting can be decomposed
into two parts: one where the dynamic consists on periodic and
almost periodic motions ($\mathcal{I},\,\mathcal{R})$ with the
diffeomorphism acting equicontinuously, and another one where the
dynamics is expansive and similar to the hyperbolic case.
Moreover, given a set having dominated decomposition, it is
characterized its stable and unstable set, its continuation by
perturbation, and their basic pieces (see \cite{PS2}). Let us say also, that in
solving the above problem, another kind of differentiable
dynamical problem arose: how is affected the dynamics of a system,
when its smoothness is improved?

Putting theorem $1$ and $2$ together, we prove the conjecture of
Palis for surface diffeomorphisms in the $  C^1-$topology:

 {\bf Theorem 4} (\cite{PS1}){\bf : } {\em Let $M^2$ be a two
dimensional compact manifold and let $f\in Diff^1(M^2).$ Then, $f$
can be $C^1$-approximated  either by a diffeomorphism exhibiting a
homoclinic tangency or by an Axiom A diffeomorphism.}

Similar arguments, prove that the variation of the topological
entropy leads to the unfolding of homoclinic tangencies. Moreover
the presence of infinitely many sinks with unbounded period also
implies the unfolding of tangencies (see \cite{PS2}, \cite{PS4}).

We want to emphasize that theorem $4$ and the previous comments are
strictly $C^1$ (on the other hand, theorem $2$, and $3$ assume that
the map is $C^2$), and nothing is known in the $ C^2-$topology. We
would like to understand what happens in the $ C^2-$topology, since
many rich dynamical phenomena take place for smooth maps.  About
this problem we would like to make some remarks.

Recall that for smooth one-dimensional endomorphisms, the absence
of critical points was enough to guarantee hyperbolicity. But for
surface maps (and due to the lack of understanding) we required
$C^1$-robust absence of tangencies. Many of the tools used in the
$C^1$ case ($C^1$-closing Lema, perturbation of the tangent map
along finite orbits) are unknown in higher topology and even in
same particular situation they are also false (see \cite{G},
\cite{PS2}). So, taking in mind the scenario for one-dimensional
dynamics, instead of ask about the $ C^r-$robust absence of
tangencies (for $r\geq 2$), we could try to know what is the two
dimensional phenomena whose {\it presence breaks the domination}
and whose {\em absence guarantee it}? In other words, what are
{\it two dimensional critical points}?

To address these problems we should consider previously some other
weaker questions. Observe that any dominated splitting is a
continuous one. Is it true the converse, at least generically?
To answer this question we would face the following:
if there is a continuous invariant splitting over the
Limit set for an smooth map (or even assuming an stronger
hypothesis: existence of two continuous invariant foliations) can
we describe the dynamic of $f$? It is clear that are dyanmics
exhibiting continuous splitting which are not dominated, for instance,
maps exhibiting some kind of saddle connection and maps
on the torus obtained as $(x,y)\to (x+\alpha, y+\beta)$.
Are those dynamics the unique ones that do
not exhibit domination?

On the other hand, splitting dealing with critical behaviors
(tangencies or "almost tangencies") are
well known in the measure-theoretical setting. This is the case of
the non-uniform hyperbolicity (or Pesin theory), where the tangent
bundle splits for points a.e. with respect to some invariant
measure, and vectors are asymptotically contracted or expanded in
a rate that may depend on the base point. Are the invariant
measures for smooth maps on surfaces with one non zero Lyapunov
exponents, non-uniformly hyperbolic?
In other words, one non-zero lyapunov exponent implies that the
other is also non-zero?  This is not true in
general, since there exist ergodic invariant measures with only
one non zero Lyapunov exponent: measure supported on invariant circle normally
hyperbolic; measure over a non-hyperbolic periodic point;
time one map of a Cherry flows; two dimensional
version of one dimension phenomena like infinitely renormalizable
and absorbing cantor sets. Are those dynamics the unique
counterexamples?

Now we would like to say a few words about the proof of Theorem
$2$. First, it is showed that the local unstable (stable) sets are
one dimensional manifolds tangent to the direction $F$ ($E$
respectively). This is achieved by, given an {\it explicit
characterization of the Lyapunov stable sets under the hypothesis
of domination} (latter we will give a precise statement). After
that, it is proved that the length of the negative iterates of the
local unstable (positive for the local stable) manifolds are
sumable, and using arguments of distortion hyperbolicity is
concluded.

We say that the point $x$ is {\em Lyapunov stable} (in the future) if given
$\e>0$ there exists $\de>0$ such that $f^n(B_\de(x))\subset
B_\e(f^n(x))$ for any positive integer $n.$ Its characterization is done in any
dimension whit the solely assumption that one of the subbundles is
one dimensional. To avoid confusion, we call such splitting {\it
codimension one dominated splitting}.

{\bf Theorem 5: }{\em Let $f:M\to M$ be a $C^2$-diffeomorphism of
a finite dimensional compact riemannian manifold $M$ and let $\La$
be a set having a codimension one dominated splitting. Then there
exists a neighborhood $V$ of $\La$ such that if $f^n(x)\in V$ for
any positive integer $n$ and $x$ is Lyapunov stable, one of the
following holds:

$1 .$ $\omega(x)$ is a periodic orbit,

$2 .$ $\omega(x)$ is a periodic curve normally attractive supporting and
irrational rotation.}

As we said before, this theorem have important consequences related
to the direction $F$.  Its local invariant tangent manifold either
is dynamically defined (it is a subset of the local unstable set)
or there are well understood phenomena: either there are periodic
curves normally contractive $\gamma$ with small length,
or there are semi-attracting periodic points, or there are closed
invariant curves normally hyperbolic with dynamics conjugated to an
irrational rotation. With this characterization in mind, and
assuming domination over the whole manifold, it is proved that $F$
is also uniquely integrable.

\section{Robust transitivity} \label{robust}

\vskip-5mm \hspace{5mm}

As we said in the beginning, in dimension higher or equal than
$3$, there exist $ C^r-$open set of diffeomorphisms which are
transitive and non-hyperbolic ($r \geq 1$). Observe that this
phenomena take place in the $C^1$-topology, fact that it is
unknown for one of the dynamics phenomena that we consider in the
previous section: the residual sets of infinitely many sinks for
surface maps.

The first examples of robust non-hyperbolic systems (examples of
robust transitive systems which are not Anosov) were given by M.
Shub (see \cite{Sh}), who considered on the $4$-torus,
skew-products of an Anosov with a Derived of Anosov
diffeomorphisms. Then, R. Ma\~n\'e (see \cite{M}) reduced the
dimension of such examples by showing that certain Derived of
Anosov diffeomorphisms on the $3$-torus are robust transitive.
Later, L. D\'iaz, (see \cite{D1}) constructed examples obtained as
a bifurcation of an heteroclinic cycle (cycle involving points of
different indices). This last ideas was pushed in \cite{BD} where
it was  showed a general geometric construction of robust
transitive attractors.

All these systems show a {\it partial hyperbolic} splitting, which
allows the tangent bundle to split into $Df$-invariant subbundles
$TM=E^s\oplus E^c\oplus E^u,$ where the behavior of vectors in
$E^s, E^u$ under iterates of the tangent map is similar to the
hyperbolic case, but vectors in $E^c$ may be neutral for the
action of the tangent map. On the other hand, recently, it was
proved by C. Bonatti and M. Viana that there are opens sets of
transitive diffeomorphisms exhibiting a dominated splitting which
do not fall into the category of partially hyperbolic ones (see
\cite{BV}).

These new situations lead to ask two natural questions: Is there a
characterization of robust transitive sets that also gives
dynamical information about them? Can we describe the dynamics
under the assumption of either partial hyperbolicity or dominated
decomposition?

The next result shows that this two questions are extremely
related. In fact, some kind of dynamics on the
tangent bundle is implied by the robust transitivity (see \cite{M}
for surfaces, \cite{DPU} for three dimensional manifolds , and
\cite{BDP} for the $n .$dimensional case):

{\bf Theorem 6:} {\em Every robustly transitive set of a
$C^1$-diffeomorphism has dominated splitting whose extremal
bundles are uniformly volume contracting or expanding.}

The central idea in this theorem is to show that, in the lack of
domination, the eigenspaces of a linear map (obtained by
multiplying many bounded linear maps) are very unstable: by small
perturbation of each of the factors, one can mix the eigenvalues
in order to get a homothety, which will correspond to the creation
of either a sink or a source, situation not allowed in the case of
the robust transitivity. This last theorem also can be formulated
in the following way:

{\bf Theorem 7} {\cite{BDP}) {\bf :} {\em There is a residual
subset of $C^1$-diffeomorphisms such that that for any
diffeomorphism in the residual set, it is verified that for any
homoclinic class of a periodic point (the closure of the
intersection of the stable and unstable manifold of it) either has
dominated splitting or it is contained in the closure of
infinitely many sources or sinks.}

What about the converse of theorem $6$? Is it true that
generically a transitive system exhibiting some kind of splitting
is robust transitive?
On the other hand, all the examples  of robust transitivity are
based in either a property of the initial system or in a
geometrical construction. But, does exist a necessary and
sufficient condition among the partial hyperbolic system such that
transitivity is equivalent to robust transitivity? {\em Can this
property be characterized in terms of the dynamic of the tangent
map?} This is clear for Anosov maps, where transitivity implies
robust transitivity, but what about for the non-hyperbolic?. In
the direction to understand this problem, in \cite{PS5} was
introduced a dynamic on the tangent bundle enough to guarantee
robustness of transitivity.

{\bf Theorem 8: }{\em Let $f\in Diff^r(M)$ be a transitive partial
hyperbolic system  verifying that there is $n_0 > 0$ such that for
any $x\in M$ there are $y^u(x)\in W^{uu}_1(x)$ and $y^s(x)\in
W^{ss}_1(x)$ with:

$1 .$$|Df^{n_0}_{|E^c(f^m(y^u(x)))}| > 2$ for any $m>0$,

$2 .$ $|Df^{-n_0}_{|E^c(f^{-m}(y^s(x)))}| > 2$ for any $m>0$,

\noindent then,  $f$ is a non-hyperbolic robust transitive system.}

Is this property generically necessary?

All these questions naturally push in the direction to understand
the dynamic induced by either a partial hyperbolic system or a
dominated splitting. In particular, do they exhibit (generically)
spectral decomposition, as was showed for a hyperbolic system and
for domination on surfaces? Observe that the non-hyperbolicity of
these systems is related with the presence of points of different
index. Is this (generically) a necessary condition for
non-hyperbolicity? In other words, assuming that there are only
points of the same index, can one conclude (generically)
hyperbolicity? For the case of domination, it is possible to show
that if the extremal directions are one-dimensional, then they
behave topologically as a hyperbolic one (see \cite{PS4}).
Moreover, it is showed that homoclinic classes with codimension
one dominated splitting and contractive bundle are generically
hyperbolic.  But, what happens if the external directions are
one-dimensional? And what about the central directions? Of course,
this question can be considered in a simpler situation: a partial
hyperbolic splitting with only one dimensional central direction.
Does the dynamic over the central direction characterize the kind
of partial hyperbolic systems?

Many of the questions done for partial hyperbolic systems can be
formulated for {\em Iterated Function Systems.} In same sense,
these systems works as a model of partial hyperbolic ones. And its
solution, could give an indication how to deal in the general
case.

In dimension higher than two, another kind of homoclinic
bifurcation breaks the hyperbolicity: the so called heteroclinic
cycles (intersection of the stable and unstable manifolds of
points of different indices, see \cite{D1} and \cite{D2}). In
particular,  the unfold of these cycles imply the existence of
striking dynamics being the more important, the appearance of
non-hyperbolic robust transitive sets (see \cite{KP} also for
superexponential growth of periodic point associated to the unfolding of
heteoclinic cycles).
Moreover, any non-hyperbolic robust transitive sets exhibits
generically heteroclinic cycles. In same sense, these cycles play
the role for the partial hyperbolic theory as transversal
intersection play for the hyperbolic theory.

A similar conjecture as the one for surfaces, was formulated by
Palis in any dimension: {\em Every $f\in Diff^r(M), r\ge 1,$ can
be $C^r$-approximated by a diffeomorphism exhibiting either a
homoclinic tangency, a heteroclinic cycle or by one which is
hyperbolic.}

A similar approach as the one done in dimension $2$ could be done:
first, try to find the dynamic on the tangent bundle for systems
$  C^1-$ far from tangencies. About this, in \cite{LW} it was proved a
similar result as the one for surfaces: {\em far from tangencies
implies domination.} Does far from heteroclinic cycles imply
hyperbolicity? Does this imply that sets with periodic points of
different index can not accumulate one on the other? And as we
asked before: sets showing dominated decomposition exhibiting
points of the same index are generically hyperbolic? These
problems are also related with the problems involving tangencies
and sinks: can a systems showing infinitely many sinks be
approximated by another one showing tangencies? It was showed that
this is true for surfaces maps (see \cite{PS3}), and in the case
of higher dimension in \cite{PS4} is given a positive answer
assuming that the sinks accumulate on a sectional dissipative
homoclinic class.

On the other hand, there is a vast works about conservative
partial hyperbolic systems describing, in same particular cases,
their ergodic properties. The description of the dynamics strongly
use the invariance of the volume measure, information that it is
not available in the general case that we would like describe. For
references about it see the complete review on this subject done
by Burns, Pugh , Shub and Wilkinson (\cite{BPSW}). Moreover it is
showed in \cite{Bo} and \cite{BoV} some kind of dichotomy (as the
one done in theorem $8$) for conservative maps in terms of
Lyapunov exponents and domination. Also I would like to mention a
recent and remarkable work of F. Rodriguez Hertz (\cite{R}) where
is proved that many Linear automorphisms on $T^4$ are stable ergodic, using
different kinds of techniques that even if only work in the
conservative case, they could be useful to understand the general
case.

\section{Flows} \label{flows}

\vskip-5mm \hspace{5mm}

For flows, a striking example is the Lorenz attractor \cite{Lo},
given by the solutions of the polynomial vector field in $R^3$:
\begin{equation}
\label{e0} X(x,y,z)= \left\{ \begin{array}{ll}
                         \dot{x} = -\alpha x + \alpha y
\\
                         \dot{y}= \beta x -y - xz
\\
                         \dot{z}=-\gamma z + xy \, ,
\end{array} \right.
\end{equation}
where $\alpha, \beta, \gamma$ are real parameters. Numerical
experiments performed by Lorenz (for $\alpha=10, \beta=28 \mbox{
and } \gamma=8/3$ ) suggested the existence, in a robust way,  of
a strange attractor toward which tends a full neighborhood of
positive trajectories of the above system. That is, the strange
attractor could not be destroyed by any perturbation of the
parameters. Most important, the attractor contains an equilibrium
point $(0,0,0)$, and periodic points accumulating on it, and hence
can not be hyperbolic. Notably, only now, three and a half decades
after this remarkable work, was it proved \cite{Tu} that the
solutions of (\ref{e0}) satisfy such a property for values
$\alpha, \beta, \gamma$ near the ones considered by Lorenz.

However, already in the mid-seventies, the existence of robust
non-hyperbolic attractors was proved for flows introduced in
\cite{ABS} and \cite{Gu}, which we now call geometric models for
Lorenz attractors. In particular, they exhibit, in a robust way,
an attracting transitive set with an equilibrium (singularity).
Moreover, the properties of this geometrical models, allow one to
extract very complete dynamical information. A natural question raises,
is such features present for any robust transitive set?

In \cite{MPP} a positive answer for this question is given:

{\bf Theorem 9 :} {\em $C^1$ robust transitive sets with
singularities on closed $3$-manifolds verifies:

$1 .$  there are either proper attractors or proper repellers;

$2 .$ the eigenvalues at the singularities satisfy the same inequalities as
the corresponding ones at the singularity in a Lorenz geometrical
model;

$3 .$ there are partially hyperbolic
with a volume expanding central direction.}

The presence of a singularity prevents these attractors from
being hyperbolic. But they exhibit a weaker form of hyperbolicity {\em
singular hyperbolic splitting}.  This class of vector fields
contains the Axiom A systems, the geometric Lorenz attractors and
the singular horseshoes in (\cite{LP}), among other systems.
Currently, there is a rather satisfactory and complete description
of singular hyperbolic vector fields defined on 3-dimensional
manifolds (but the panorama in higher dimensions remains open).
More precisely, it is proved in a sequel of works that a singular
hyperbolic set for flow is  $K^* .$expansive, the periodic orbits
are  dense in its limit set, and it has a spectral decomposition
(see \cite{PP}, \cite{K}).

On the other hand, for the case of flows, appears a new kind of bifurcation that leads to a new dynamics distinct
from the ones for diffeomorphism: the so called {\em singular cycles} (cycles involving singularities and periodic
orbits, see \cite{BLMP}, \cite{Mo}, \cite{MP} and \cite{MPP1} for examples of dynamics in the sequel of the
unfolding of it). Systems exhibiting this cycles are dense among open set of systems exhibiting a singular
hyperbolic splitting. Moreover, recently A. Arroyo and F. Rodriguez Hertz (see \cite{AR}) studying the dynamical
consequences of the dominated splitting for the Linear Poincare flow,  proved that {\em any three dimensional flow
can be $  C^1-$approximated either by a flow exhibiting tangency or singular cycle, or by a hyperbolic one.}

\label{lastpage}

\end{document}